\documentclass [a4paper,leqno,12pt]{article}
\usepackage[intlimits]{amsmath}
\usepackage{amsbsy}
\usepackage{amstext}
\usepackage{amssymb}
\usepackage[english]{babel}
\usepackage{indentfirst}
\usepackage{a4wide}
\newtheorem{theorem}{Theorem}
\newtheorem{proposition}[theorem]{Proposition}
\newtheorem{lemma}[theorem]{Lemma}

\newtheorem{cor}[theorem]{Corollary}
\numberwithin{theorem}{section} \numberwithin{equation}{section}

\def\proof{{\bf Proof. }}
\def\endprf {\mbox{$\blacksquare$}}
\newcommand{\beq}{\begin{equation}}
\newcommand{\eeq}{\end{equation}}

\DeclareMathOperator*{\argmax}{argmax}

\DeclareMathOperator*{\meas}{meas}
\DeclareMathOperator*{\essinf}{essinf}

\title{Optimal relocation strategies for spatially mobile
consumers} \date{2nd May 2007}
\author{Iordan Iordanov\footnote{Faculty of Mathematics and Informatics, Sofia University, e-mail iordanov@fmi.uni-sofia.bg} \and Andrey Vassilev\footnote{Faculty of Mathematics and Informatics, Sofia University, e-mail avassilev@fmi.uni-sofia.bg}}

\begin{document} \maketitle \begin{abstract}
We develop a model of the behaviour of a dynamically optimizing
economic agent who makes consumption-saving and spatial relocation
decisions. We formulate an existence result for the model, derive
the necessary conditions for optimality and study the behaviour of
the economic agent, focusing on the case of a wage distribution
with a single maximum.
\\ {}  \end{abstract}
\section{Introduction}\label{sec:intro}

The emergence of the literature on ``new economic geography'' in
the 1990s has rekindled the interest in the spatial aspects of
economics. The new generation of models makes heavy use of the
standard economics toolkit and analyzes a number of issues from a
dynamic perspective or from the perspective of optimizing agents.
Interestingly, however, spatial models adopting the perspective of
\textit{dynamically optimizing} consumers remain in relative
minority, despite the fact that they are standard fare in
mainstream economic research. The models in \cite{Bal01},
\cite{Bou06} and \cite{Bri04} are notable exceptions in this
respect.

The present work develops a model that studies the behaviour of a
dynamically optimizing economic agent who makes two types of
interrelated choices: consumption-saving decisions and spatial
relocation (migration) decisions. Unlike the constructs in
\cite{Bal01}, \cite{Bou06} and \cite{Bri04}, the consumer in our
model has a finite lifetime and a bequest incentive at the end of
his life. This departure from classical Ramsey-type models enables
richer global dynamics by allowing agents to inherit their
ancestors' savings in a setup akin to that of overlapping
generations models. It also offers the additional option of
introducing heterogeneous agents whose economy-wide behaviour can
be obtained through an explicit aggregation rule.

A second important difference with the above papers is that our
consumer saves in nominal assets. This partly depends on our
choice to center the model around the behaviour of (potentially
different) individuals as opposed to that of a representative
agent. More importantly, however, the nominal savings feature
reflects our belief that pecuniary considerations play an
important role in the choice of where to work and how much to
consume.

Formally, we cast the model in the form of a continuous-time
optimal control problem with a finite planning horizon. The
assumptions of the model, while fairly standard in economics,
create several mathematical challenges in the present setup.
First, they preclude the direct application of the existence
theorems for optimal control problems known to the authors. This
requires an alternative approach to proving the existence of
solutions of the model. In particular, unlike traditional
existence proofs in the spirit of Theorem 4, \S 4.2 in
\cite{Lee69}, we prove an existence result that dispenses with
convexity assumptions on the set of generalized speeds for the
optimal control problem. Also, the functional forms employed in
the model do not allow one to directly apply Pontryagin's maximum
principle, since the transversality condition for one of the state
variables is not defined at the point $0$. To be able to use the
maximum principle, we prove that for an optimal control-trajectory
pair the terminal value of the particular state variable is
strictly positive. Finally, economic considerations point to the
fact that only a subset of the possible values for the other state
variable in the model are of real interest. One way to take care
of that issue is to constrain the values of this state variable to
lie in a certain set -- the interval $[0,1]$ in our case -- for
each point in time. However, instead of using an explicit state
constraint, which would complicate the use of the maximum
principle, we introduce an additional correcting mechanism by
suitably defining the wage distribution function $w(x)$ outside
the interval $[0,1]$. We claim this mechanism does not influence
the other characteristics of the model, while being sufficient to
ensure that the optimal state variable never leaves the set in
question, and we prove that indeed this is the case.

The rest of the paper is organized as follows. Section
\ref{sec:model} introduces the model and the assumptions we make.
Section \ref{sec:exist} proves the existence of a solution to the
model under the above assumptions. Section \ref{sec:NC} applies
Pontryagin's maximum principle to obtain necessary conditions for
optimality. Section \ref{sec:NCtransf} describes some convenient
transformations of the system of necessary conditions and comments
on the existence of solutions to this system. The analysis in
section \ref{sec:termasset} characterizes the asymptotic behaviour
of terminal assets for different sets of model parameters. This
establishes facts that are useful for the study of relocation
choices in section \ref{sec:reloc}. The results in this section
are also of independent economic interest as they shed light on
the impact of intra- and intertemporal preferences on the saving
decisions of an individual with a sufficiently long planning
horizon. Finally, section \ref{sec:reloc} tackles the question of
relocation behaviour in the basic case of a wage distribution
having a single maximum (single-peaked wage distribution). The
results obtained for this case are intuitive, if unsurprising: in
most cases a consumer with a sufficiently long lifespan relocates
toward the wage maximum. While the single-peak case offers easily
predictable results, we consider it useful as a testing ground for
the model before applying it to more interesting situations.
Indeed, preliminary results by the authors on the case of a
double-peaked wage distribution suggest that a host of complex
situations, including multiple solutions and bifurcations, can
arise.
\\
\\
\textbf{Acknowledgements.} We would like to thank Tsvetomir
Tsachev and Vladimir Veliov for useful discussions and pointers to
the literature. The responsibility for any errors is solely ours.

\section{The model}\label{sec:model}

We employ a continuous-time model that deals with the case of a
consumer who, given an initial location in space $x_0$ and asset
level $a_0$, supplies inelastically a unit of labour in exchange
for a location-dependent wage $w(x(t))$, and chooses consumption
$c(t)$ and spatial location $x(t)$ over time. The consumer has a
finite lifetime $T$ at the end of which a bequest in the form of
assets is left. This bequest provides utility to the consumer.
More precisely, for  $\rho$, $r$, $\eta$, $\xi$ and $p$ --
positive constants, and $\theta \in (0,1)$, we look at the optimal
control problem \beq \max_{c(t),z(t) \in \Delta} J(c(t),z(t)) :=
\int_{0}^{T}e^{-\rho t}\left(
\frac{c(t)^{1-\theta}}{1-\theta}-\eta z^2(t) \right)dt+e^{-\rho
T}\frac{a(T)^{1-\theta}}{1-\theta} \label{eq:obj} \eeq subject to
\beq \dot{a}(t)=ra(t)+w(x(t))-p c(t) -\xi z^2(t),
\label{eq:assets} \eeq \beq
\dot{x}(t)=z(t),\label{eq:location}\eeq
$$a(0)=a_0\geq 0,$$
$$x(0)=x_0\in [0,1],$$
where $a(t)$, $x(t)$ are the state variables, assumed to be
absolutely continuous, and $c(t)$, $z(t)$ are the control
variables. The set of admissible controls $\Delta$ consists of all
pairs of functions $(c(t),z(t))$ which are measurable in $[0,T]$
and satisfy the conditions \beq 0 \leq c(t) \leq C,
\label{eq:cconstr}\eeq \beq |z(t)|\leq Z, \label{eq:zconstr} \eeq
\beq a(T)\geq 0. \label{eq:aTnonneg}\eeq The constants $C$ and $Z$
are such that \beq C^{\theta} >
\max\left(1,\mu^{\frac{1}{\theta}}\right) \left( \frac{a_0 +
T\max_x |w(x)|}{p\frac{1-e^{-rT}}{r}}\right) ,
\label{eq:Ccapbound} \eeq \beq Z
> \frac{T \max_x |w'(x)| e^{rT}}{2\xi},\label{eq:Zcapbound}\eeq where $\mu := \max_{t,t_0 \in
[0,T]}e^{(r-\rho)(t-t_0)}>0$.

\textbf{Remark.} The bounds we impose on the admissible controls
through equations \eqref{eq:cconstr} and \eqref{eq:zconstr} are
convenient from a technical viewpoint when proving the existence
theorem in section \ref{sec:exist}. Conditions
\eqref{eq:Ccapbound} and \eqref{eq:Zcapbound} ensure that these
constraints are never binding. However, considerations of general
nature -- both economic and physical -- make such constraints
appealing.

In the above model $\rho>0$ is a time discount parameter and
$\theta \in (0,1) $ is the utility function parameter. The control
$c(t)$ represents physical units of consumption and the control
$z(t)$ governs the speed of relocation in space. We assume that
relocation in space brings about two type of consequences. First,
relocation causes subjective disutility associated with the fact
that there is habit formation with respect to the place one
occupies. Second, changing one's location is associated with
monetary relocation costs that have to be paid out of one's income
or stock of assets. As a baseline case we choose to capture these
phenomena by means of the speed of movement in space $\dot{x}(t)$
or, equivalently, $z(t)$, transformed through a quadratic
function. The manner in which spatial relocation affects the
consumer's utility and wealth can vary widely, however, therefore
other functional forms are certainly admissible. The parameters
$\eta,\xi \geq 0$ multiplying this function measure the subjective
disutility from changing one's location in space and the
relocation costs in monetary terms, respectively. The parameters
$p>0$ and $r>0$ stand for the price of a unit of consumption and
the interest rate, respectively.

The nonnegativity condition is imposed on terminal assets $a(T)$
both to have a well-defined objective functional and to capture
the intuitive observation that, with a known lifetime, a debtor is
unlikely to be allowed to leave behind outstanding liabilities to
creditors. The condition $a(T)\geq 0$ also sheds light on the
nonnegativity restriction for $a_0$, since in an environment where
no debts are allowed at the end of one's lifetime, no debtor
position can be inherited at birth.

For the purposes of our analysis we look at the basic case where
economic space is represented by the real line. We are interested
in only a subset of it, the interval $[0,1]$. This is modelled by
taking the initial location $x_0 \in [0,1]$ and requiring the
location-dependent wage, which is positive in $(0,1)$, to be
negative outside $[0,1]$ and to satisfy additional assumptions.
Namely, we have $w(x)> 0,~x \in (0,1)$ and $w(x)<0,~x \not\in
[0,1]$, as well as $w'(x)>0$, $x \in (-\infty,0]$ and $w'(x)<0$,
$x\in [1,\infty)$. Later in the paper we formally verify the
intuitive claim that an optimal trajectory for $x(t)$ will never
leave the interval $[0,1]$ under the above conditions. We also
assume that $w(x) \in C^2(\mathbb{R}^1)$ and $w(x)$ is bounded,
i.e. $\max_{x \in \mathbb{R}}|w(x)|<+\infty$. We impose additional
requirements on $w(x)$ to derive some of the results in section
\ref{sec:reloc}.

\section{Existence of solutions}\label{sec:exist}

Next we investigate the issue of existence of a solution to the
model. The proof requires two intermediate results, shown as
lemmas below.

\begin{lemma} Let the functions $x_i,~i=1,2,\ldots,$ and $\bar{x}$ be defined on $[0,T]$ and take values in the interval $[a,b]$. Let $x_i$ tend uniformly to $\bar{x}$ as $i \rightarrow \infty$ (denoted by $x_i \rightrightarrows \bar{x}$) and $w \in C^0[a,b]$. Then, in $[0,T]$, as $m\rightarrow \infty$ we have
\begin{enumerate}
  \item[i)] $\frac{1}{m}\sum_{i=1}^m x_i \rightrightarrows
  \bar{x}$,
  \item[ii)] $w(x_m)\rightrightarrows w(\bar{x})$,
  \item[iii)] $w(\frac{1}{m}\sum_{i=1}^m x_i)\rightrightarrows
  w(\bar{x})$.
\end{enumerate}
\label{thm:lem1}\end{lemma}

\proof The proof directly replicates the standard proofs of
counterpart results on numerical sequences.
\endprf

\begin{lemma}[The Banach-Saks Theorem] Let $\{ v_n \}_{n=1}^{\infty}$ be a sequence of elements in a Hilbert space $H$ which are bounded in norm: $\| v_n \|\leq K = const,~\forall n \in \mathbb{N}$. Then, there exist a subsequence $\{ v_{n_k} \}_{k=1}^{\infty}$ and an element $v\in H$ such that $$\left\| \frac{v_{n_1}+\cdots+v_{n_s}}{s}-v \right\|\rightarrow 0 \textrm{ as } s\rightarrow \infty .$$ \label{thm:lem2}\end{lemma}
\proof See, for example, \cite[pp.78-81]{Die75}. \endprf

\begin{theorem} Under the assumptions stated in section \ref{sec:model}, there exists a solution $(c(t),z(t))\in \Delta$ of problem \eqref{eq:obj}-\eqref{eq:location}. \label{thm:exist}\end{theorem}
\proof We start by noting that the set of admissible controls
$\Delta$ is nonempty. To see this, choose controls $c(t)\equiv
c_0=const$ and $z(t)\equiv 0$. Then, any $c_0\in (0, w(x_0)/p]$
will ensure that $a(T)\geq 0$.

Next, observe that $(c(t),z(t))\in \Delta$ implies $c(t),z(t)\in
L_\infty[0,T]$ and \beq 0 \leq a(T) \leq const = e^{rT}(a_0+T
\max_{x \in \mathbb{R}}|w(x)|). \label{eq:aTbounded}\eeq (We note
that \eqref{eq:aTnonneg} implies the following bounds,  $$ p \|
c(t) \|_{L_1[0,T]},\xi \| z(t) \|^2_{L_2[0,T]} \leq e^{rT}\left(
a_0+T \max_x |w(x)| \right) , $$which do not depend on the
constants $C$ and $Z$.)

Through an application of H\"{o}lder's inequality one verifies
that $\int_0^T c(t)^{1-\theta}e^{-\rho t}dt \leq const(T) \|
c(t)\|^{1-\theta}_{L_1} $.

Thus, for $(c(t),z(t))\in \Delta$, the objective functional
\eqref{eq:obj} is bounded. Consequently,
$J_0:=\sup_{(c(t),z(t))\in \Delta}J(c(t),z(t))<\infty$. Then we
can choose a sequence of controls $\{ (c_k(t),z_k(t)) \}\subset
\Delta$ such that $J(c_k(t),z_k(t))\rightarrow J_0$.

Let $a_k(t)$ and $x_k(t)$ be the state variables corresponding to
the controls $(c_k(t),z_k(t))$. It is easy to verify that the
functions $a_k(t)$ and $x_k(t)$ form a uniformly bounded and
equicontinuous set. Then, by the Arzel\`a-Ascoli theorem (see,
e.g., \cite{Lee69}, Ch.4), there exists a subsequence
$(a_{k_s}(t),x_{k_s}(t))\rightrightarrows
(\bar{a}(t),\bar{x}(t))$.

Then, if $c_{k_s}(t)$ and $z_{k_s}(t)$ are the controls
corresponding to $(a_{k_s}(t),x_{k_s}(t))$, by Lemma
\ref{thm:lem2} we can in turn choose subsequences $c_{k_{s_q}}(t)$
and $z_{k_{s_q}}(t)$ whose arithmetic means tend in $L_2[0,T]$
norm to some elements in $L_2[0,T]$, denoted $\bar{c}(t)$ and
$\bar{z}(t)$, respectively. However, we do not claim that
$\bar{a}(t)$ and $\bar{x}(t)$ correspond to $\bar{c}(t)$ and
$\bar{z}(t)$. For brevity we introduce the notation
$c_q(t):=c_{k_{s_q}}(t)$, $z_q(t):=z_{k_{s_q}}(t)$ etc., as well
as $\tilde{c}_m(t):=\frac{1}{m}\sum_{q=1}^{m}c_q(t)$ and
$\tilde{z}_m(t):=\frac{1}{m}\sum_{q=1}^{m}z_q(t)$.

Then, we have established that: (1)
$(a_q(t),x_q(t))\rightrightarrows (\bar{a}(t),\bar{x}(t))$ as $q
\rightarrow \infty$ and (2)
$\tilde{c}_m(t)\xrightarrow[L_2]{}\bar{c}(t)$,
$\tilde{z}_m(t)\xrightarrow[L_2]{}\bar{z}(t)$ as $m\rightarrow
\infty$.

Recall that $a_q(t)$ and $x_q(t)$ correspond to $c_q(t)$ and
$z_q(t)$ as solutions to the respective differential equations
\eqref{eq:assets} and \eqref{eq:location}.

So far, it is not clear whether $\tilde{c}_m(t)$ and
$\tilde{z}_m(t)$ are admissible. It is immediately seen that they
satisfy \eqref{eq:cconstr} and \eqref{eq:zconstr} but the
corresponding $a(T)$ may fail to satisfy \eqref{eq:aTnonneg}.
However, we can show that the controls $\bar{c}(t)$ and
$\bar{z}(t)$ are admissible.

To prove the last claim, note first that according to \cite[Ch.7,
\S 2.5, Prop.4]{Kol76} we can choose a subsequence of
$\{\tilde{c}_m(t), \tilde{z}_m(t)\}$ that converges a.e. to
$(\bar{c}(t),\bar{z}(t))$ and, after passing to the limit, we
obtain that $\bar{c}(t)$ and $\bar{z}(t)$ satisfy
\eqref{eq:cconstr} and \eqref{eq:zconstr}.

It remains to show that $\bar{\bar{a}}(T)=e^{rT}\left[a_0+\int_0^T
[w(\bar{x}(t))-p \bar{c}(t)-\xi \bar{z}^2(t)]e^{-rt}dt \right]\geq
0,$ where $\bar{\bar{x}}(t)=x_0+\int_0^t \bar{z}(\tau) d\tau$.

Consider \beq \tilde{a}_m(T)=e^{rT}\left[a_0+\int_0^T
[w(\tilde{x}_m(t))-p \tilde{c}_m(t)-\xi\tilde{z}^2_m(t)]e^{-rt}dt
\right],\label{eq:a_mT}\eeq with $\tilde{x}_m(t)=x_0+\int_0^t
\tilde{z}_m(\tau) d\tau=\frac{1}{m}\sum_{q=1}^m \left(
x_0+\int_0^t z_q(\tau) d\tau \right)=\frac{1}{m}\sum_{q=1}^m
x_q(t)$. Adding and subtracting $\frac{1}{m}\sum_{q=1}^m
w(x_q(t))$, and applying Jensen's inequality to the term
$\tilde{z}^2_m(t)$, we obtain \beq
\begin{split} \tilde{a}_m(T) \geq & e^{rT}\int_0^T\left[
w(\tilde{x}_m(t))- \frac{1}{m}\sum_{q=1}^m
w(x_q(t))\right]e^{-rt}dt + \\ & \frac{1}{m}\sum_{q=1}^m
e^{rT}\left[ a_0+\int_0^T [w(x_q(t))-p c_q(t)-\xi
z^2_q(t)]e^{-rt}dt \right] \geq \\ & e^{rT}\int_0^T\left[
w(\tilde{x}_m(t))- \frac{1}{m}\sum_{q=1}^m
w(x_q(t))\right]e^{-rt}dt
\end{split}\label{eq:a_mT1}
\eeq By Lemma \ref{thm:lem1} both integrands inside the square
brackets in the last line of \eqref{eq:a_mT1} tend uniformly to
$w(\bar{x}(t))$, so that the integral tends to zero. Thus, if
$\lim_{m\rightarrow \infty}\tilde{a}_m(T)$ exists, we have
$\lim_{m\rightarrow \infty}\tilde{a}_m(T) \geq 0$.

We proceed to check that $\lim_{{m_j}\rightarrow
\infty}\tilde{a}_{m_j}(T)=\bar{\bar{a}}(T)$ for a suitable
subsequence $\tilde{a}_{m_j}(T)$. We know that
$\frac{1}{m}\sum_{q=1}^m
x_q(t)=x_0+\int_0^t\frac{1}{m}\sum_{q=1}^m z_q(\tau)d \tau$. Since
$\frac{1}{m}\sum_{q=1}^m x_q(t) \rightrightarrows \bar{x}(t)$ and,
additionally, it is easy to verify by applying H\"{o}lder's
inequality that $\int_0^t\tilde{z}_m(\tau)d \tau \rightarrow
\int_0^t \bar{\bar{z}}(\tau)d\tau$ when
$\tilde{z}_m(t)\xrightarrow[L_2]{}\bar{z}(t)$, we obtain
$\bar{x}(t)=x_0+\int_0^t \bar{z}(\tau)d\tau =\bar{\bar{x}}(t)$.

As $\tilde{c}_m(t)\xrightarrow[L_2]{}\bar{c}(t)$ and
$\tilde{z}_m(t)\xrightarrow[L_2]{}\bar{z}(t)$, there exist a
subsequences $\tilde{c}_{m_j}(t)$ and $\tilde{z}_{m_j}(t)$ such
that $\tilde{c}_{m_j}(t)\xrightarrow[a.e.]{}\bar{c}(t)$ and
$\tilde{z}_{m_j}(t)\xrightarrow[a.e.]{}\bar{z}(t)$. To simplify
notation, we refer to the new subsequences as $\tilde{c}_j(t)$ and
$\tilde{z}_j(t)$. Since the function $z^2$ is bounded on $[-Z,Z]$,
by Lebesgue's dominated convergence theorem $\int_0^T \xi
\tilde{z}^2_{j}(t)e^{-rt}dt\rightarrow \int_0^T \xi
\bar{z}^2(t)e^{-rt}dt$. It can also be verified that $\int_0^T
\tilde{c}_{j}(t)e^{-rt}dt\rightarrow \int_0^T
\bar{c}(t)e^{-rt}dt$. Lastly, we know that $\int_0^T
w(\tilde{x}_{j}(t))e^{-rt}dt \rightarrow \int_0^T
w(\bar{x}(t))e^{-rt}dt$ as $w(\tilde{x}_{j}(t)) \rightrightarrows
w(\bar{x}(t))$. Consequently, the limit of \eqref{eq:a_mT} as $m_j
\rightarrow \infty$ exists and is equal to $\bar{\bar{a}}(T)$, so
that $\bar{\bar{a}}(T) \geq 0$. This shows that $\bar{c}(t)$ and
$\bar{z}(t)$ are admissible.

By an application of Lebesgue's dominated convergence theorem to
the respective terms in \eqref{eq:obj}, we get $\lim_{j
\rightarrow
\infty}J(\tilde{c}_j(t),\tilde{z}_j(t))=J(\bar{c}(t),\bar{z}(t))$.

Define $\rho_{m_j}(T) := e^{rT}\int_0^T\left[
w(\tilde{x}_{m_j}(t))- \frac{1}{m_j}\sum_{q=1}^{m_j}
w(x_q(t))\right]e^{-rt}dt$. Obviously,
$\tilde{\tilde{a}}_{m_j}(T)=\tilde{a}_{m_j}(T)-\rho_{m_j}(T)$ also
tends to $\bar{\bar{a}}(T)$ and $\tilde{\tilde{a}}_{m_j}(T) \geq
\frac{1}{m_j}\sum_{q=1}^{m_j} a_q(T)$, where $a_q(T)$ corresponds
to $(c_q(t),z_q(t))$. Then, indexing by $j$ instead of $m_j$ to
simplify notation, we get

$$J_0 \geq
J(\bar{c}(t),\bar{z}(t))=\lim_{j \rightarrow \infty}\left\{
\int_0^T \left[\frac{\tilde{c}_j(t)^{1-\theta}}{1-\theta}-\eta
\tilde{z}^2_j(t)\right]e^{-\rho t }dt + e^{-\rho
T}\frac{\tilde{\tilde{a}}^{1-\theta}_j(T)}{1-\theta} \right\} \geq
$$

$$\lim_{j \rightarrow \infty}\left\{\frac{1}{j}
\sum_{i=1}^{j}\left[\int_0^T
\left[\frac{c^{1-\theta}_i(t)}{1-\theta}-\eta
z^2_i(t)\right]e^{-\rho t }dt + e^{-\rho T}\frac{
a^{1-\theta}_i(T)}{1-\theta}\right]\right\}=$$

$$ \lim_{j \rightarrow \infty}\left\{\frac{1}{j}
\sum_{i=1}^{j}J(c_i(t),z_i(t)) \right\} = J_0,$$ where the
inequality is a consequence of the fact that the functions $\sigma
\mapsto \sigma^{1-\theta}$ and $z \mapsto (-z^2)$ are concave and
we can apply Jensen's inequality. This shows that the admissible
pair $(\bar{c}(t),\bar{z}(t))$ is optimal, as required.
\endprf

\section{Necessary conditions for optimality}\label{sec:NC}

In this section we turn to the derivation of a set of necessary
conditions for optimality on the basis of Pontryagin's maximum
principle. To apply the maximum principle, however, we need to
ensure that the terminal utility from assets $e^{-\rho
T}a(T)^{1-\theta}/(1-\theta)$ is well-behaved at least for the
optimal value of terminal assets. To this end, we prove the
following

\begin{theorem} For the optimal controls $(c(t),z(t))$ the terminal value of
assets $a(T)$ is strictly positive for any $T>0$.
\label{thm:aT>0}\end{theorem} \proof Let us assume that there is a
time $T_0>0$ for which $a(T_0)=0$.

\textbf{Step 1.} We first verify that it is impossible to have
$c(t) \equiv 0$. Assuming that $c(t) \equiv 0$, together with
$a(T_0)=0$, yields the objective functional
$$J(0,z(t))=-\eta \int_0^{T_0}z^2(t)e^{-\rho t}dt \leq 0.$$

If one of the following two conditions is valid:
\begin{enumerate}
\item $a_0>0$ and $x_0 \in [0,1]$;
\item $a_0=0$ and $x_0 \in (0,1)$,
\end{enumerate} then we can choose the admissible pair $\bar{z}(t)\equiv
0$ (so that $x(t)\equiv x_0$) and $\bar{c}(t)\equiv c_0 = const
>0$, where $c_0$ is such that $$a_0+\int_0^{T_0}[w(x_0)-p c_0]e^{-rt}dt=0.$$ The last
condition is equivalent to $$a_0+T_0 w(x_0)\frac{e^{-r
T_0}-1}{-r}=p c_0 \frac{e^{-r T_0}-1}{-r}$$ and therefore $c_0>0$.
Then
$$J(\bar{c}(t),\bar{z}(t))=\int_0^{T_0}\frac{c_0^{1-\theta}}{1-\theta}e^{-\rho t}dt
>0,$$ contradicting the optimality of $(c(t),z(t))$.

The case $a_0=0$ and $x_0=0$ or $1$ is pathological in the sense
that the consumer has neither current income (w(0)=w(1)=0), nor
initial wealth. Economically, it is implausible to expect that
such a consumer will manage to obtain a loan. From a purely formal
point of view, however, the consumer could get a loan and finance
his relocation even in this case. Moreover, he will be able to
attain positive consumption levels.

To illustrate the above claim, suppose that $x_0=0$, $a_0=0$ and
the consumer spends all the income left after paying the
relocation costs. Fix $\varepsilon_0 > 0$ in such a way that
$w'(x) \geq w'(0)/2 >0$ for $x \in [0,\varepsilon_0]$. Let the
relocation strategy be given by the control
$\bar{z}(t)=\varepsilon \sin \frac{\pi}{T}t,~\varepsilon>0$. Then
consumption is given by $\bar{c}(t)=w(\bar{x}(t))-\xi
\bar{z}^2(t)$, where $\bar{x}(t)$ is
$$\bar{x}(t)=\varepsilon \int_0^t \sin \left( \frac{\pi}{T}\tau \right)d\tau=\frac{\varepsilon T}{\pi}\left( 1-\cos \frac{\pi t}{T} \right)=\frac{2 \varepsilon T}{\pi}\sin^2 \frac{\pi t}{2 T}. $$
Then, $\bar{x}(T)=\frac{2T \varepsilon}{\pi}<\varepsilon_0$ for
$\varepsilon$ sufficiently small. Notice that $$w(\bar{x}(t))
=w(\bar{x}(t))-w(0)=w'(x^*(t))\bar{x}(t)\geq
\frac{w'(0)}{2}\bar{x}(t),$$ for some $x^*(t)\in (0,\bar{x}(t))$.
Consequently, we obtain
$$w(\bar{x}(t))-\xi \bar{z}^2(t)\geq \varepsilon \left[ \frac{w'(0)}{2}\frac{2 T}{\pi}\sin^2 \frac{\pi t}{2 T}-\varepsilon \xi \sin^2 \frac{\pi t}{T}  \right] = \varepsilon \sin^2 \frac{\pi t}{2 T}\left[ \frac{T w'(0)}{\pi}- 4\varepsilon \xi \cos^2 \frac{\pi t}{2T}  \right]. $$
Consumption will be positive if $$g(t):=\frac{T w'(0)}{\pi}-
4\varepsilon \xi \cos^2 \frac{\pi t}{2T}>0 \textrm{ for } t \in
[0,T].$$ For $\varepsilon$ small $g(0)=T w'(0)/\pi -4\xi
\varepsilon > 0$. Also, $$g'(t)=4\varepsilon \xi \frac{\pi}{2T}2
\cos \frac{\pi t}{2T} \sin \frac{\pi t}{2T} = 2 \varepsilon \xi
\frac{\pi}{T} \sin \frac{\pi t}{T} \geq 0 \textrm{ for }t \in
[0,T].$$ Thus, $g(t) \geq g(0)>0$, as required.

\textbf{Remark.} It is easy to see that in the above example we
can take $\bar{z}(t)$ to be any smooth function that is positive
on $(0,T)$, zero for $t=0,T$ and $\dot{\bar{z}}(0)>0$.

\textbf{Step 2.} Since $c(t) \not \equiv 0$, there exists a set $A
\subset [0,T] $, $\meas A>0$, such that $$\essinf_{t \in A} c(t)>
\varepsilon_1>0.$$

Let us take the control pair $(\bar{c}(t),\bar{z}(t))$ with
$\bar{c}(t):= c(t)-\varepsilon \chi_A(t)$ and $\bar{z}(t) :=
z(t)$, where $\chi_A(t)$ is the indicator function of the set $A$
and $\varepsilon \in (0,\varepsilon_1)$. These controls are
admissible if we have terminal assets $\bar{a}(T_0)>0,~\forall
\varepsilon \in(0,\varepsilon_1)$. To verify the last claim, we
take
$$\bar{a}(T_0)=e^{r T_0}\left[
a_0+\int_0^{T_0}[w(x(s))-p(c(s)-\varepsilon \chi_A(s))-\xi
z^2(s)]e^{-rs}ds \right]=e^{r T_0}\int_A p \varepsilon
e^{-rs}ds=\varepsilon C_1,$$ where $C_1:=p e^{r T_0} \int_A
e^{-rs}ds >0$.

An application of Taylor's formula yields
$$\frac{\bar{c}(t)^{1-\theta}}{1-\theta}=\frac{c(t)^{1-\theta}}{1-\theta}+(-\varepsilon \chi_A(t))c(t)^{-\theta}+(-\varepsilon \chi_A(t))^2
\frac{-\theta}{2}c^*(t)^{-\theta-1},$$ where $c^*(t)=\alpha(t)
\bar{c}(t)+(1-\alpha(t))c(t),~\alpha(t) \in (0,1)$ or
$c^*(t)=c(t)-\varepsilon \chi_A(t) \alpha(t)$. Note also that for
$t \in A$ we have $0<c(t)-\varepsilon_1\cdot 1 \leq c^*(t)\leq
c(t)$, so that $(c(t)-\varepsilon_1)^{-\theta-1}\geq
c^*(t)^{-\theta-1}\geq c(t)^{-\theta-1}$.

Let us compare
$$J(c(t),z(t))=\int_0^{T_0}\frac{c(t)^{1-\theta}}{1-\theta}e^{-\rho t}dt-\eta\int_0^{T_0}z^2(t)e^{-\rho
t}dt$$ and
\begin{equation*} \begin{split} J(\bar{c}(t),\bar{z}(t))= & \int_0^{T_0}\frac{\bar{c}(t)^{1-\theta}}{1-\theta}e^{-\rho
t}dt-\eta\int_0^{T_0}z^2(t)e^{-\rho
t}dt+\frac{\bar{a}(T_0)^{1-\theta}}{1-\theta}e^{-\rho T_0} \\ = &
J(c(t),z(t)) + \left[ -\varepsilon \int_A c(t)^{-\theta}e^{-\rho
t}dt -\frac{\theta \varepsilon^2}{2}\int_A(c(t)-\varepsilon
\alpha(t))^{-1-\theta}e^{-\rho t}dt  \right]+ \\ & +
\frac{(\varepsilon C_1)^{1-\theta}}{1-\theta}e^{-\rho T_0}.
\end{split}
\end{equation*}

We will show that $J(\bar{c}(t),\bar{z}(t))>J(c(t),z(t))$ for
$\varepsilon \in (0,\varepsilon_1)$ sufficiently small. This will
be true if we are able to establish that $$\frac{(\varepsilon
C_1)^{1-\theta}}{1-\theta}e^{-\rho T_0} > \varepsilon \int_A
c(t)^{-\theta}e^{-\rho t}dt + \frac{\varepsilon^2 \theta}{2}\int_A
(c(t)-\varepsilon_1 )^{-1-\theta}e^{-\rho t}dt ,$$ where the last
integral provides an upper bound on $\int_A (c(t)-\varepsilon
\alpha(t))^{-1-\theta}e^{-\rho t}dt $. Denoting the respective
positive constants in the above inequality by $B_1$, $B_2$ and
$B_3$, we obtain
$$\varepsilon^{1-\theta}B_1 > \varepsilon B_2 + \varepsilon^2 B_3
$$ or $$B_1 > \varepsilon^\theta B_2 + \varepsilon^{1+\theta}
B_3,$$ which is obviously true for $\varepsilon \in
(0,\varepsilon_1)$ sufficiently small. This contradicts the
optimality of $(c(t),z(t))$. Thus, $a(T_0)=0$ cannot be true and
hence $a(T_0)>0$.
\endprf

On the basis of Theorem \ref{thm:aT>0} an optimal solution
$\bar{c}(t),\bar{z}(t)$ to problem
\eqref{eq:obj}-\eqref{eq:location} (possibly non-unique) also
solves the following problem, where the controls $(c(t),z(t))\in
\Delta_1 \subset \Delta$:

$$\max_{c(t),z(t) \in \Delta_1}
J(c(t),z(t))$$
$$\dot{a}(t)=ra(t)+w(x(t))-p c(t) -\xi z^2(t)$$
$$\dot{x}(t)=z(t)$$
$$a(0)=a_0\geq 0,$$
$$x(0)=x_0\in [0,1],$$
$$a(T)\geq \delta >0,$$ with $\delta$ being an appropriate
constant, strictly smaller than the optimal value of terminal
assets.

To avoid burdensome notation, from now on we do not append
additional symbols to the state, costate and control variables in
the model when referring to their optimal values. However, we use
alternative symbols to denote alternative sets of variables to be
compared with the optimal ones.

Taking into account that we do not impose any state constraints on
the problem, Theorem 5.2.1 in \cite{Cla83} provides the set of
necessary conditions. To derive the latter, we define the
Hamiltonian for the problem \beq
\begin{split}H & := H(t,a,x,\varphi,\psi,p_1,p_2)= \\  & p_1(r
a+w(x)-p \varphi -\xi \psi^2 )+p_2 \psi+\lambda_0 e^{-\rho
t}\left( \frac{\varphi^{1-\theta}}{1-\theta} -\eta \psi^2
\right),\end{split} \label{eq:Hamilt}\eeq where $\lambda_0 \in
\{0,1\}$. Then

1) The costate variables $p_i(t),~i=1,2,$ satisfy a.e. on $(0,T)$
\beq \dot{p}_1(t)=-r p_1(t), \label{eq:p1dot}\eeq \beq
\dot{p}_2(t)=-p_1(t)w'(x(t)). \label{eq:p2dot}\eeq

2) The function $\varphi,\psi \mapsto
H(t,a(t),x(t),\varphi,\psi,p_1(t),p_2(t))$ attains its maximum
with respect to $\varphi,\psi$ at the point $(c(t),z(t))$ for
almost all $t \in [0,T]$, where $\varphi,\psi$ satisfy the
constraints on the function values arising from $\Delta_1$, i.e.
$\varphi \in [0,C]$ and $|\psi|\leq Z$: \beq
\begin{split}
H(t) & := H(t,a(t),x(t),c(t),z(t),p_1(t),p_2(t)) =  \\
& \max_{\varphi, \psi } H(t,a(t),x(t),\varphi,\psi,p_1(t),p_2(t)).
\end{split} \label{eq:Hammax}\eeq

3a) Since $a(t)$ and $x(t)$ are fixed at $t=0$, the values of
$p_i(0)$ are arbitrary, i.e. \beq p_1(0)=\lambda_1,~ \lambda_1 \in
\mathbb{R} , \label{eq:trans3}\eeq  \beq p_2(0)=\lambda_2,~
\lambda_2 \in \mathbb{R} . \label{eq:trans4}\eeq

3b) Since the terminal values $(a(T),x(T))$ of the state variables
are at an interior point of the target set $$\{ (a,x)\in
\mathbb{R}^2 | a \geq \delta >0, x \in \mathbb{R}^1 \},$$ the
corresponding normal cone is trivial and the transversality
condition at the right endpoint $T$ has the form \beq
p_1(T)=\lambda_0 e^{-\rho T}a(T)^{-\theta}, \label{eq:trans1}\eeq
\beq p_2(T)=0 , \label{eq:trans2}\eeq (cf. condition 4) in Theorem
5.2.1  and the functional form for $f(x(b))$ in \S 5.2 in
\cite{Cla83}).

4) The variables $p_1(t),p_2(t),\lambda_0$ are not simultaneously
equal to zero.

Below we specify the form of the necessary conditions in greater
detail.

According to \eqref{eq:p1dot} and \eqref{eq:trans3} we have \beq
p_1(t)=\lambda_1 e^{-rt}. \label{eq:p1solved}\eeq

\begin{proposition}If there exists $t_0 \in [0,T]$ such that $c(t_0)\in
(0,C)$, then $c(t)>0$ for almost all $t\in [0,T]$.
\label{thm:dtpositive}
\end{proposition}
\proof If there exists $t_0$ with the above properties, then
\eqref{eq:Hammax} implies
$$\frac{\partial}{\partial\varphi}H(t_0,a(t_0),x(t_0),\varphi,z(t_0),p_1(t_0),p_2(t_0))\left |_{\varphi=c(t_0)}=0,
\right.$$ i.e. \beq -p \lambda_1 e^{-r t_0}+\lambda_0 e^{-\rho
t_0}c(t_0)^{-\theta} = 0. \label{eq:supplHammax1}\eeq

If we assume that $\lambda_0=0$, then $\lambda_1=0$ and, because
of \eqref{eq:p1solved}, one obtains $p_1(t)\equiv 0$. This implies
that $p_2(t)\equiv const =\lambda_2$. Now \eqref{eq:trans2} shows
that $\lambda_2=0$, which constitutes a contradiction with
condition 4) from the cited theorem in \cite{Cla83}. Therefore,
$\lambda_0=1$.

Assume that there exists $t_1 \in [0,T]$ for which $c(t_1)=0$.
Then, for all sufficiently small $\varphi >0$ we have
$$\frac{H(t_1,a(t_1),x(t_1),\varphi,z(t_1),p_1(t_1),p_2(t_1))-H(t_1,a(t_1),x(t_1),0,z(t_1),p_1(t_1),p_2(t_1))}{\varphi-0}\leq
0,$$ i.e. $$-p_1 (t_1) p + \lambda_0 e^{-\rho t_1}
\frac{\varphi^{-\theta}}{1-\theta}\leq 0.$$ Since $\lambda_0>0$,
for $\varphi \rightarrow 0+$ the last inequality leads to a
contradiction. This proves the proposition.
\endprf

\begin{cor} If there exists $t_1 \in [0,T]$ for which $c(t_1)=0$, then $\lambda_0=0$ and $c(t)=0$ for almost all $t \in [0,T]$. \label{thm:lambda0zero}\end{cor}
\proof The conclusion on $\lambda_0$ can be obtained  in the same
manner as in the proof of Proposition \ref{thm:dtpositive} by
passing to the limit as $\varphi \rightarrow 0+$. If we assume the
existence of a point $t_0\in [0,T]$ for which $c(t_0)>0$, we can
proceed as in the proof of the proposition and get $p_1(t)\equiv
p_2(t)\equiv 0$ and $\lambda_0=0$, which is impossible. \endprf

\begin{proposition}The optimal consumption cannot be identically
zero. \label{thm:conspositive}
\end{proposition}
\proof Assume that the controls $c(t)\equiv 0$ and $z(t)$ are
optimal. Then $$J(0,z(t))=-\int_0^T \eta z^2(t)e^{-\rho t}dt+
e^{-\rho T}\frac{a(T)^{1-\theta}}{1-\theta}.$$ Take the controls
$\tilde{c}(t)=\varepsilon$ and $\tilde{z}(t)=z(t)$, where
$\varepsilon>0$ is sufficiently small. These controls are
admissible, as the respective value of terminal assets is
$$\tilde{a}(T)=e^{rT}\left\{ a_0+\int_0^T [w(x(t))-p\varepsilon -\xi z^2(t)]e^{-rt}dt  \right\}=a(T)-\varepsilon
C_1,$$where $C_1:=e^{rT}p\int_0^Te^{-rt}dt>0$. It is evident that
for $\varepsilon$ sufficiently small we have
$\tilde{a}(T)>\delta$, since $a(T)>\delta$. It remains to check
that for $\varepsilon$ close to zero we have
\begin{equation*}\begin{split} J(\varepsilon,z(t))= &  \int_0^T
\frac{\varepsilon^{1-\theta}}{1-\theta}e^{-\rho t}dt-\eta \int_0^T
z^2(t)e^{-\rho t}dt+ e^{-\rho T}\frac{(a(T)-\varepsilon
C_1)^{1-\theta}}{1-\theta}
> \\ & J(0,z(t))=-\int_0^T \eta z^2(t)e^{-\rho t}dt+
e^{-\rho T}\frac{a(T)^{1-\theta}}{1-\theta}, \end{split}
\end{equation*} which is equivalent to $$\varepsilon^{1-\theta} C_2 > \frac{e^{-\rho T}}{1-\theta}\left[a(T)^{1-\theta}-(a(T)-\varepsilon
C_1)^{1-\theta}\right],~C_2:=const>0.$$ The last expression is
obviously true for all $\varepsilon$ sufficiently small.
\endprf

\textbf{Remark.} So far it is clear  that the optimal consumption
cannot be identically zero and that if there exists $t_0$ such
that $c(t_0)\in (0,C)$, then $\lambda_0=1$. It remains to check
whether we can have $c(t)=C$ for some $t$.

We first establish the following result. \begin{proposition} It is
impossible for the optimal $c(t)$ to satisfy \beq c(t)\geq C_0>0,
\label{eq:dtgeqD1}\eeq \label{thm:dtgeqD1} where \beq C_0 >
\frac{a_0 + T\max_x |w(x)|}{p\frac{1-e^{-rT}}{r}} .
\label{eq:D1def}\eeq
\end{proposition}
\proof Notice that if condition \eqref{eq:dtgeqD1} is true, then
the inequality $a(T)\geq \delta$ is violated. Indeed, if
\eqref{eq:dtgeqD1} holds, then
\begin{equation*}\begin{split} a(T) = & e^{rT}\left[ a_0+\int_0^T
\left[w(x(t))-pc(t) -\xi z^2(t)\right]e^{-rt}dt \right]\leq \\ &
e^{rT}\left[ a_0+\int_0^T \left[\max_x |w(x)|-p C_0
\right]e^{-rt}dt \right],\end{split}\end{equation*} which is
negative when \eqref{eq:D1def} holds. \endprf

\begin{proposition} The number $\lambda_1$ is strictly positive.
\label{thm:lambda1positive}\end{proposition} \proof We know that
for the optimal $c(t)$ it is impossible to have $c(t)\geq C_0$ or
$c(t)\equiv 0$. Consequently, there exists $t_0\in [0,T]$ for
which $c(t_0) \in (0,C_0)$. Then
$$\frac{\partial}{\partial\varphi}H(t_0,a(t_0),x(t_0),\varphi,z(t_0),p_1(t_0),p_2(t_0))\left
|_{\varphi=c(t_0)}=0, \right.$$ and hence \eqref{eq:supplHammax1}
holds. This in turn implies that $\lambda_0=1$, as well as $$p
\lambda_1 e^{-r t_0}  =  e^{-\rho t_0}c(t_0)^{-\theta}.$$
Therefore, we have $\lambda_1>0$ and \beq \lambda_1=
\frac{e^{(r-\rho)t_0}c(t_0)^{-\theta}}{p}.
\label{eq:lambda1expr}\eeq
\endprf

\begin{proposition} There does not exist $t\in [0,T]$ for which $c(t)=C$.
\label{thm:dtlessthanD}\end{proposition}\proof Assuming the
contrary, by the maximum principle we obtain
$$\frac{H(t,a(t),x(t),\varphi,z(t),p_1(t),p_2(t))-H(t,a(t),x(t),C,z(t),p_1(t),p_2(t))}{\varphi-C}\geq
0$$ for $\varphi \in (0,C)$ and so for $\varphi\rightarrow C-0$ we
get $$- p\lambda_1 e^{-rt} + e^{-\rho t} C^{-\theta}\geq 0,$$
which implies
\begin{equation*}
C^{\theta}\leq \frac{e^{(r-\rho)t}}{\lambda_1
p}=\frac{e^{(r-\rho)t}c(t_0)^{\theta}}{e^{(r-\rho)t_0}}<
e^{(r-\rho)(t-t_0)}C_0^\theta \leq \mu C_0^\theta,
\end{equation*}
where $\mu := \max_{t,t_0 \in [0,T]}e^{(r-\rho)(t-t_0)}>0$. In
other words, $$C \leq \mu^{\frac{1}{\theta}}C_0,$$ which is
impossible. \endprf

The results obtained so far allow us to to find an expression for
the optimal consumption $c(t)$.
\begin{cor}For each $t \in [0,T]$ we have $c(t)\in(0,C)$. The optimal consumption rule has the form \beq c(t)=\left[
\frac{1}{p \lambda_1
}\right]^{\frac{1}{\theta}}e^{\frac{r-\rho}{\theta} t
}=\frac{1}{p^{\frac{1}{\theta}}}e^{\frac{\rho-r}{\theta}(T-t)}a(T).
\label{eq:copt1}\eeq \label{thm:copt} \end{cor}

Before deriving an expression for the optimal relocation control
$z(t)$, we note that \eqref{eq:p2dot} and \eqref{eq:trans2} imply
\beq p_2(t)=\lambda_1 \int_t^T w'(x(\tau))e^{-r\tau}d\tau =
e^{(r-\rho)T}a(T)^{-\theta}\int_t^T w'(x(\tau))e^{-r\tau}d\tau.
\label{eq:p2trans1} \eeq

\begin{proposition} For each $t \in [0,T]$ we have the strict inequality $$|z(t)|< Z.$$ \label{thm:zbound} \end{proposition}
\proof Assume, for example, that $z(t_1)=Z$ for some $t_1 \in
[0,T]$. Then, after passing to the limit in the respective
difference quotient, we obtain
$$-p_1(t_1)\xi 2Z+p_2(t_1)-2\eta Z e^{-\rho t_1}\geq
0,$$ so that $$p_2(t_1)\geq 2(\xi\lambda_1 e^{-rt_1}+\eta e^{-\rho
t_1 })Z.$$

Similarly, the assumption that $z(t_2)=-Z$ for some $t_2 \in
[0,T]$ leads to $$-p_2(t_2)\geq 2(\xi \lambda_1 e^{-r t_2}+\eta
e^{-\rho t_2 })Z.$$

In both cases we have ($i=1$ or $2$)
\begin{equation*}
\begin{split}
Z  \leq & \frac{\pm p_2(t_i)}{2(\xi \lambda_1 e^{-r t_i}+\eta
e^{-\rho t_i })}  \leq  \frac{|p_2(t_i)|}{2(\xi \lambda_1 e^{-r
t_i}+\eta e^{-\rho t_i })}  \leq  \frac{ \lambda_1
\left|\int_{t_i}^T w'(x(\tau))e^{-r\tau}d\tau \right| }{2(\xi
\lambda_1 e^{-r t_i}+\eta e^{-\rho t_i })}  \leq \\ & \\ & \frac{
\max_x \left| w'(x) \right| |T-t_i| }{2(\xi e^{-r
t_i}+\frac{\eta}{\lambda_1} e^{-\rho t_i })} < \frac{T \max_x
|w'(x)|}{2\xi e^{-r t_i}} \leq \frac{T \max_x |w'(x)|}{2\xi} e^{r
T},
\end{split}
\end{equation*} which is impossible by the definition of $Z$.
\endprf

\begin{cor}
For $t\in [0,T]$ we have for the optimal relocation speed $z(t)\in
(-Z,Z)$ and then \beq z(t)=\frac{p_2(t)}{2 (\xi \lambda_1 e^{-rt}+
\eta e^{-\rho t})}. \label{eq:zopt1} \eeq\label{thm:zopt}
\end{cor}

\section{Existence of a solution of the system of necessary
conditions}\label{sec:NCtransf}

In order to facilitate the study of the differential equations
arising from the set of necessary conditions in section
\ref{sec:NC}, it would prove convenient to rewrite the
differential system. Theorem \ref{thm:exist} guarantees the
existence of a solution to the problem
\eqref{eq:obj}-\eqref{eq:location} which in turn ensures the
existence for each $T>0$ of a solution to the  following problem:
\beq \left |
\begin{array}{l} \dot{x}(t)=\frac{y(t)}{F(t)}, \\ \dot{y}(t)=-w'(x(t))\lambda_1 e^{-rt}, \\ x(0)=x_0, \\ y(T)=0, \end{array}
\right. \label{eq:AuxSys1}\eeq where $y(t):=p_2(t)$ and
$F(t):=2(\xi \lambda_1 e^{-rt}+\eta e^{-\rho t})$. It follows that
there exists a solution to the problem \beq \left |
\begin{array}{l} \frac{d}{dt}\left( F(t) \dot{x}(t)  \right)+w'(x(t))\lambda_1 e^{-rt}=0, \\ x(0)=x_0, \\ \dot{x}(T)=0. \end{array}
\right. \label{eq:AuxSys2}\eeq

The latter fact can also be established without recourse to
Theorem \ref{thm:exist}. Following the procedure described in \S
73 of \cite{Lov24}, we construct the respective Green function and
transform \eqref{eq:AuxSys2} in the form \beq x(t)=\int_0^T K(t,
\tau)\lambda_1 e^{-r \tau}w'(x(\tau))d\tau, \label{eq:Lovitt}\eeq
where $$K(t, \tau)=\left \{ \begin{array}{l l}
\int_0^\tau \frac{1}{F(s)}ds, & \tau \in [0,t] \\
\int_0^t \frac{1}{F(s)}ds, & \tau \in [t,T]
\end{array} \right .$$

Since the function $w'(x)$ is bounded and continuous by
assumption, a solution to \eqref{eq:Lovitt} exists. This is a
consequence of Leray-Schauder index theory (see \S 2.4 in
\cite{Nir74}). Also, the solution to \eqref{eq:AuxSys2} may not be
unique, as can be seen from simple examples of eigenfunction
problems that possess nontrivial solutions.

A solution to \eqref{eq:AuxSys1} or \eqref{eq:AuxSys2} can be
viewed as a particular member of the family of solutions
$(x(t,\alpha),y(t,\alpha))$ to the Cauchy problem \beq \left |
\begin{array}{l} \dot{x}(t)=\frac{y(t)}{F(t)}, \\ \dot{y}(t)=-w'(x(t))\lambda_1 e^{-rt}, \\ x(0)=x_0, \\ y(0)=\alpha, \end{array}
\right. \label{eq:AuxSys3}\eeq where $\alpha$ has been chosen
appropriately, so that \beq y(T,\alpha)=0. \label{eq:AuxSys4}\eeq
The existence of a unique solution to \eqref{eq:AuxSys3} on the
interval $[0,T]$ for initial data $(x_0,\alpha)$ and each $T>0$ is
ensured by Corollary 3.1, chapter 2, in \cite{Har64}.

Since \eqref{eq:AuxSys4} is equivalent to $\dot{x}(T,\alpha)=0$,
we can integrate the differential equation in \eqref{eq:AuxSys2}
over $[0,T]$ to arrive at an equivalent form of
\eqref{eq:AuxSys4}: \beq \alpha = \lambda_1 \int_0^T w'(x(\tau,
\alpha))e^{-r\tau}d\tau . \label{eq:AuxSys5}\eeq

It is straightforward to verify the following
\begin{proposition}The function $x(t)\equiv x_0$ is a solution to \eqref{eq:AuxSys2} if and only if the point $x_0$ is a critical point for $w(x)$, i.e. $w'(x_0)=0$.
\label{thm:5.1}\end{proposition} The analysis of the solutions of
the system of necessary conditions, which is carried out in
section \ref{sec:reloc}, provides the dynamics of the behaviour of
the economic agent, implied by this model, in the baseline case
when the wage distribution has a single maximum point on the
interval $[0,1]$. Prior to that, the next section studies the
properties of the function $T \mapsto a(T)$ as $T \rightarrow
\infty$.

\section{Dynamics of terminal assets $a(T)$ for different time horizons}\label{sec:termasset}

In this section we study the dependence of optimal terminal assets
$a(T)$ on the length of the time horizon $T$. Although terminal
assets is the natural object of study due to the fact that it is
easily interpretable in economic terms, the discussion may equally
well be framed in terms of the behaviour of $\lambda_1$, viewed as
a function of the time horizon $T$ and denoted $\lambda_1(T)$.
This approach is feasible by virtue of the relationship \beq
\lambda_1(T)=e^{(r-\rho)T}a(T)^{-\theta}.\label{eq:lambda1T}\eeq
Below we derive upper and lower bounds on $\lambda_1(T)$, which
will be needed in the analysis of section \ref{sec:reloc}. We
assume for simplicity that $x_0 \in (0,1)$, i.e. $w(x_0)>0$, as
well as that $a_0>0$. Also, in this section we denote by $C$
different constants that do not depend on $T$. Since we do not use
the bound on the control $c(t)$ from \eqref{eq:cconstr} in this
section, no confusion can arise from this convention.

\subsection{An upper bound on $\lambda_1(T)$}\label{sec:termassup}

The pair $(c(t)\equiv c_0=const,z(t)\equiv 0)$ is admissible for
$c_0$ appropriately chosen. Then $x(t) \equiv x_0$ and we set
$c_0:=\frac{w(x_0)}{p}$. In this case terminal assets are
$$\tilde{a}(T)=e^{rT}\left[
a_0+\int_0^T[w(x(t))-pc(t)-\xi  z^2(t)]e^{-rT}dt \right]=a_0
e^{rT}.$$ Consequently, \begin{equation*} J(c_0,0) =
c_1(1-e^{-\rho T})+c_2 e^{[r(1-\theta)-\rho]T},
\end{equation*}where $c_1$ and $c_2$ are constants that depend on
$c_0$ and $a_0$.

On the other hand, for the optimal controls $(c(t),z(t))$ we have
\begin{equation*}
J(c(t),z(t)) =  \int_0^T
\left[a(T)\frac{1}{p^{1/\theta}}e^{\frac{\rho-r}{\theta}(T-t)}
\right]^{1-\theta}\frac{e^{-\rho t}}{1-\theta} dt-\eta\int_0^T
z^2(t)e^{-\rho t}dt+\frac{a(T)^{1-\theta}}{1-\theta}e^{-\rho T},
\end{equation*} which takes different forms depending on whether
$\rho-r(1-\theta)$ is different from zero.

Since $J(c(t),z(t))\geq J(c_0,0)$, for $\rho-r(1-\theta)\neq 0$ we
obtain $$\frac{a(T)^{1-\theta}}{1-\theta}e^{-\rho T}\left[
1+\frac{1}{p^{\frac{1-\theta}{\theta}}}\frac{e^{\frac{\rho-r(1-\theta)}{\theta}T}-1}{\frac{\rho-r(1-\theta)}{\theta}}
\right]\geq c_1+c_2e^{-(\rho-r(1-\theta))T}-c_1 e^{-\rho T}.$$

Thus, \beq a(T) \geq \left\{ \frac{\tilde{c}_1 e^{\rho
T}+\tilde{c}_2 e^{r(1-\theta)T}-\tilde{c}_1
}{1+\frac{1}{p^{\frac{1-\theta}{\theta}}}\frac{e^{\frac{\rho-r(1-\theta)}{\theta}T}-1}{\frac{\rho-r(1-\theta)}{\theta}}}
\right\}^{\frac{1}{1-\theta}}, \label{eq:aTlowbound}\eeq where
$\tilde{c}_1=(1-\theta)c_1$, $\tilde{c}_2=(1-\theta)c_2$. Using
the last expression together with \eqref{eq:lambda1T}, we can
derive upper bounds on $\lambda_1(T)$.

\begin{proposition}
Under the assumptions of this section, we have ($\forall T > 0$):
\beq \lambda_1(T) \leq \left\{
  \begin{array}{l l}
    C, & \textrm{if } \rho-r(1-\theta)>0, \\
    C e^{-(\rho-r(1-\theta))T}, & \textrm{if } \rho-r(1-\theta)<0, \\
    C (1+T)^{\frac{\theta}{1-\theta}}, & \textrm{if }
    \rho-r(1-\theta)=0.
  \end{array}
 \right.
\label{eq:lambda1Tuprbnd}\eeq and, accordingly, \beq a(T) \geq
\left\{
  \begin{array}{l l}
    C e^{\frac{r-\rho}{\theta}T}, & \textrm{if } \rho-r(1-\theta)>0, \\
    C e^{rT}, & \textrm{if } \rho-r(1-\theta)<0, \\
    C \frac{e^{\frac{\rho}{1-\theta}T}}{(1+T)^{\frac{1}{1-\theta}}}, & \textrm{if }
    \rho-r(1-\theta)=0.
  \end{array}
 \right.
\label{eq:aTuprbnd}\eeq \label{thm:uprbndsummary}
\end{proposition}

\subsection{A lower bound on $\lambda_1(T)$}\label{sec:termasslo}

We first look at a particular case of the main problem, for which
\beq w(x)\equiv W =const. \label{eq:WageConst}\eeq Then
$\dot{p}_2\equiv 0$ which, together with the transversality
condition $p_2(T)=0$, yields $p_2(t)\equiv 0$, i.e. $z(t)\equiv 0$
and $x(t)\equiv x_0.$

The optimal consumption rule is
$c(t)=\frac{1}{p^{1/\theta}}\bar{a}(T)e^{\frac{\rho-r}{\theta}(T-t)},$
where $\bar{a}(T)$ is the optimal terminal value of assets for the
problem with condition \eqref{eq:WageConst}. We will calculate
$\bar{a}(T)$ from $$\bar{a}(T)=e^{rT}\left[
a_0+\int_0^T\left(W-p^{\frac{\theta-1}{\theta}}\bar{a}(T)e^{\frac{\rho-r}{\theta}(T-t)}\right)e^{-rt}dt
\right].$$ Thus, we find \beq \bar{a}(T) \leq \left\{
  \begin{array}{l l}
    C e^{\frac{r-\rho}{\theta}T}, & \textrm{if } \rho-r(1-\theta)>0, \\
    C e^{rT}, & \textrm{if } \rho-r(1-\theta)<0, \\
    C \frac{e^{rT}}{1+T}, & \textrm{if }
    \rho-r(1-\theta)=0.
  \end{array}
 \right.
\label{eq:abarTlwrbnd}\eeq

\begin{proposition} Let $w(x)\leq W,~\forall x,$ and let $a(T)$ and
$\bar{a}(T)$ be the optimal terminal asset values for the problems
with wage distributions $w(x)$ and $W$, respectively (all other
parameters of the two problems being identical). Then
$$a(T) \leq \bar{a}(T).$$
\label{thm:aT_leq_abarT}\end{proposition} \proof Since according
to \eqref{eq:copt1} optimal consumption for the two problems has
the form $a(T)\Psi(t)$ and $\bar{a}(T)\Psi(t)$ with one and the
same function $\Psi(t)$, we obtain
$$[\bar{a}(T)-a(T)]\left(1+e^{rT}\int_0^Tp\Psi(t)e^{-rt}dt
\right)=e^{rT}\int_0^T [W-w(x(t))]e^{-rt}dt+e^{rT}\xi \int_0^T
z^2(t)e^{-rt}dt,$$ where $x(t)$ and $z(t)$ refer to the variables
in the problem with wage distribution $w(x)$. This completes the
proof. \endprf

From Proposition \ref{thm:aT_leq_abarT} and equations
\eqref{eq:lambda1T} and \eqref{eq:abarTlwrbnd}, we obtain
\begin{proposition}
Under the assumptions of this section, we have ($\forall T > 0$):
\beq a(T) \leq \left\{
  \begin{array}{l l}
    C e^{\frac{r-\rho}{\theta}T}, & \textrm{if } \rho-r(1-\theta)>0, \\
    C e^{rT}, & \textrm{if } \rho-r(1-\theta)<0, \\
    C \frac{e^{rT}}{1+T}, & \textrm{if }
    \rho-r(1-\theta)=0.
  \end{array}
 \right.
\label{eq:aTlwrbnd}\eeq and, accordingly, \beq \lambda_1(T) \geq
\left\{
  \begin{array}{l l}
    C, & \textrm{if } \rho-r(1-\theta)>0, \\
    C e^{-(\rho-r(1-\theta))T}, & \textrm{if } \rho-r(1-\theta)<0, \\
    C (1+T)^{\theta}, & \textrm{if }
    \rho-r(1-\theta)=0.
  \end{array}
 \right.
\label{eq:lambda1Tlwrbnd}\eeq \label{thm:lwrbndsummary}
\end{proposition}

\textbf{Remark.} The bounds derived above can be refined in some
cases. For instance, the first inequality in \eqref{eq:aTlwrbnd}
implies very different behaviour of $a(T)$ depending on whether
$\rho \in (r(1-\theta),r)$, $\rho=r$ or $\rho>r$.

\section{Optimal relocation strategies for single-peaked wage distributions}\label{sec:reloc}

This section studies the optimal relocation behaviour of the
consumer, as described by $x(t)$, in the important case of
single-peaked wage distributions on the interval $[0,1]$. We
demonstrate first the validity of the following general claim
(under the conditions stated at the end of section
\ref{sec:model}):
\begin{proposition} The optimal trajectory $x(t)$ remains in the
interval $[0,1]$, regardless of the particular form of the wage
distribution $w(x)$ in $[0,1]$. \label{thm:xin01}\end{proposition}
\proof Notice that since $\dot{x}(t)=p_2(t)/F(t)$, with $F(t)$
defined as in section \ref{sec:NCtransf}, in view of
\eqref{eq:p2trans1} we can write
$$\dot{x}(t)=\frac{\lambda_1}{F(t)}\int_t^T w'(x(\tau))e^{-r\tau}d\tau=G(t) \int_t^T
w'(x(\tau))e^{-r\tau}d\tau ,$$ where $G(t):=\lambda_1/F(t)$.

Assume first that at time $t_1$ the point $x(t)$ leaves the
interval $[0,1]$ to the left (i.e. leaves the interval at $x=0$)
and remains to the left of zero until $t=T$, so that $x(t)<0$ for
$t \in (t_1,T]$. Then, for $t \in [t_1,T]$, $w'(x(t))>0$ and
consequently $\dot{x}(t)>0$. This would imply that for some $t_*
\in (t_1,T)$, $x(T)-x(t_1)=(T-t_1)\dot{x}(t_*)>0$, or
$x(T)>x(t_1)=0$, which contradicts the assumption that $x(t)<0$
for $t \in (t_1,T]$. Thus, $x(t)$ cannot leave the interval
$[0,1]$ to the left and remain outside it until the end of the
planning horizon $T$. A similar argument shows that it is
impossible for $x(t)$ to leave the interval $[0,1]$ to the right
and stay there.

Let us now assume that $x(t)$ leaves the interval $[0,1]$ to the
left of zero at time $t_1$ and returns back at time $t_2>t_1$.
Again, for $t\in (t_1,t_2)$ we have $x(t)<0$ and $w'(x(t))>0$,
which means that $\int_{t_1}^{t_2}w'(x(t))e^{-rt}dt > 0$. Since
$x(t)$ leaves the interval $[0,1]$ to the left at $t_1$, it must
be that $\dot{x}(t_1)\leq 0$. By the same logic, at time $t_2$ we
should have $\dot{x}(t_2) \geq 0$. Then one obtains
\begin{equation*}
\begin{split} \dot{x}(t_1) = & G(t_1)\int_{t_1}^{T}w'(x(t))e^{-rt}dt =
G(t_1)\int_{t_1}^{t_2}w'(x(t))e^{-rt}dt+\frac{G(t_1)}{G(t_2)}G(t_2)\int_{t_2}^{T}w'(x(t))e^{-rt}dt
\\ = &
G(t_1)\int_{t_1}^{t_2}w'(x(t))e^{-rt}dt+\frac{G(t_1)}{G(t_2)}\dot{x}(t_2)>0,
\end{split}
\end{equation*} which contradicts the condition $\dot{x}(t_1)\leq 0$.
Hence it is impossible for $x(t)$ to temporarily leave the
interval $[0,1]$ to the left. Naturally, this argument can be
applied with obvious modifications to the hypothesis that $x(t)$
temporarily goes to the right of $x=1$. Summarizing the above
conclusions, we see that the optimal $x(t)$ remains in the
interval $[0,1]$. \endprf

We turn next to the main object of study for this section: the
case when the wage function has a single peak on the interval
$[0,1]$. The example of the quadratic function $w(x)=x(1-x)$ in a
neighbourhood of the interval $[0,1]$ may facilitate
visualization.

Assume that in this case the initial location $x_0$ lies to the
left of the wage peak, i.e. if $x_1:=\argmax_{x \in [0,1]}w(x)$,
then $0 \leq x_0 < x_1$. For the remainder of this section we will
assume that $w'(x)>0,~x\in [0,x_1]$ and $w'(x)<0,~x\in [x_1,1]$.
In this case, for $T$ sufficiently small, $x(T)$ will remain in a
small neighbourhood of $x_0$. However, this means that
$w'(x(t))>0$ for any $t \in [0,T]$ and therefore $\dot{p}_2(t)<0$,
which implies $p_2(t)>0$ since $p_2(T)=0$. As $p_2(t)>0$, we
obtain $\dot{x}(t)>0$. In words, for sufficiently small planning
horizons the consumer unambiguously relocates toward the wage
maximum.

If $x(T) \in (x_0,x_1)$, we have, using the notation in section
\ref{sec:NCtransf}, $\dot{y}(t)=-w'(x(t))p_1(t)<0$. Since
$y(T)=0$, it follows that $y(t)>0,~t \in [0,T)$. Then
$\dot{x}(t)=\frac{y(t)}{F(t)}>0$, so that $x(t)$ is monotonically
increasing.

We note that there does not exist a solution to the system of
necessary conditions for which $x(T)=x_1$. For such a solution we
would have $y(T)=0$ and one could compare this solution of the
stationary solution $\tilde{x}(t) \equiv x_1, \tilde{y}(t)\equiv
0$. Then, the uniqueness of the solution to a Cauchy problem (for
identical data at $t=T$) shows that the two solutions coincide.
This, however, is impossible, since for $t=0$ the values of the
two solutions are different ($x(0)=x_0 \neq \tilde{x}(0)=x_1$).

We would like to check whether it is possible for the terminal
location $x(T)$ to lie to the right of the wage peak for
$x_0<x_1$. To this end, assume that $x(T)>x_1$ and let $t_1$ be
the time when point $x_1$ is reached last, i.e. $x(t_1)=x_1$ and
for $t \in (t_1,T)$ we have $x(t)>x_1$. (In other words, $t_1=\sup
\{ t\in [0,T]| x(t)=x_1\} $.) Then, by the mean value theorem,
$0<x(T)-x(t_1)=(T-t_1)\frac{y(t_*)}{F(t_*)}$ for some $t_* \in
(t_1,T)$. However, since $\dot{y}(t)>0$, $t\in (t_1,T)$, and
$y(T)=0$ imply $y(t_*)<0$, we obtain $\frac{y(t_*)}{F(t_*)}<0$,
which is a contradiction. Thus, $x(T)$ cannot lie to the right of
$x_1$.

The systematic study of the relocation behaviour of the economic
agent in this case can be reduced to the analysis of the way in
which the solutions to the Cauchy problem \eqref{eq:AuxSys3}
behave for different values of the parameter $\alpha$. Those
solutions that satisfy \eqref{eq:AuxSys4} are also solutions to
the system of necessary conditions \eqref{eq:AuxSys1}, i.e.
extremals. Through this approach we can also obtain information on
the number of solution to the problem at hand. Of course, if only
one extremal exists, then it is the solution we seek.

We remind the reader that the number $\lambda_1=\lambda_1(T)$ is
fixed, insofar as $T$ is fixed.

\textbf{Case I: $\alpha \leq 0$.} It is obvious that for small $t$
we have $y(t)<0$ since $\dot{y}(t)<0$. For such $t$ we have
$$x(t)=x_0+\int_0^t \frac{y(\tau)}{F(\tau)}d\tau < x_0,$$ i.e. the
agent shifts toward $x=0$. It means that $\dot{y}(t)$ remains
negative, so that $y(t)=\alpha+\int_0^t \dot{y}(\tau)d\tau$ also
remains negative and $x(t)$ keeps moving to the left. Thus, it is
impossible for \eqref{eq:AuxSys4} to become true, i.e. there are
no extremals among the solutions of \eqref{eq:AuxSys3} for $\alpha
\leq 0$.

\textbf{Case II: $\alpha > 0$.} In this case we have $y(t)>0$ in a
neighbourhood of $t=0$ and so $x(t)$ moves to the right in the
direction of the point $x_1$. However,
$\dot{y}(t)=-w'(x(t))p_1(t)<0$, so that $y(t)$ decreases. If it
turns out that $y(T)=0$ and $x(T)<x_1$, then the respective
solution is an extremal. The existence of such an extremal is
guaranteed by Theorem \ref{thm:exist}. The latter claim can be
established through an alternative approach, which allows us to
ascertain the number of extremals.

We introduce the notation $M(\alpha)$ for the right-hand side of
\eqref{eq:AuxSys5}. Since $x(t)<x_1$ for $t\in [0,T]$, we get
$$M(\alpha)\leq \lambda_1 \max_x |w'(x)| \frac{1-e^{-rT}}{r}=:M_0,$$
i.e. in view of \eqref{eq:AuxSys5} the relevant values of $\alpha$
lie in the interval $(0,M_0)$.

It is easy to see that the function $$g(\alpha) :=
\alpha-M(\alpha)$$ is continuous on the interval $[0,M_0+1]$ and
satisfies the inequalities $$g(0)<0<g(M_0+1).$$ Consequently,
there exists $\alpha>0$ for which $g(\alpha)=0$, i.e. which
satisfies \eqref{eq:AuxSys5}.

\begin{proposition} For the case of a single-peaked wage distribution $w(x)$ with $w''(x) \leq 0$ in $[0,1]$, there exists a unique extremal for the system \eqref{eq:AuxSys1}.
\label{thm:SingleExtremal}\end{proposition} \proof Assume that at
least two different extremals exist. They solve the system
\eqref{eq:AuxSys3} for different positive values $\alpha_1 \neq
\alpha_2$, for which $\alpha_i-M(\alpha_i)=0,~i=1,2$. Then \beq
(\alpha_2-\alpha_1)\left(1-\frac{d}{d\alpha}M(\alpha^*)\right)=0,~\alpha^*=\varkappa
\alpha_1+(1-\varkappa)\alpha_2,~\varkappa\in(0,1).\label{eq:diffalpha1alpha2}\eeq
The derivative
$$\left. \frac{d}{d\alpha}\left( \lambda_1 \int_0^T w'(x(t,
\alpha))e^{-rt}dt \right)\right|_{\alpha=\alpha^*} $$ has the form
$$\lambda_1 \int_0^T w''(x(t,\alpha^*))\frac{\partial
x(t,\alpha^*)}{\partial \alpha}e^{-rt}dt ,$$ with
$x_\alpha(t,\alpha) := \frac{\partial x(t,\alpha)}{\partial
\alpha} $ and $y_\alpha(t,\alpha) := \frac{\partial
y(t,\alpha)}{\partial \alpha} $ satisfying the equations of
variation \cite[Ch.V, Theorem 3.1]{Har64}:
\begin{equation*} \left |
\begin{array}{l} \dot{x}_\alpha(t,\alpha)=\frac{y_\alpha(t,\alpha)}{F(t)}, \\ \dot{y}_\alpha(t,\alpha)=-w''(x(t,\alpha))x_\alpha(t,\alpha) \lambda_1 e^{-rt}, \\ x_\alpha(0,\alpha)=0, \\ y_\alpha(0,\alpha)=1. \end{array}
\right. \end{equation*} Consequently, $x_\alpha (t,\alpha^*)$
solves the linear equation
$$\frac{d}{dt}\left(F(t)\dot{x}_\alpha(t,\alpha^*)\right)+w''(x (t,\alpha^*))\lambda_1 e^{-rt}x_\alpha
(t,\alpha^*)=0$$ for initial data $x_\alpha (0,\alpha^*)=0$ and
$\dot{x}_\alpha (0,\alpha^*)=1$. Multiplying by $x_\alpha
(t,\alpha^*)$ and integrating over $(0,t)$, we obtain
$$F(t)\dot{x}_\alpha (t,\alpha^*)x_\alpha
(t,\alpha^*) = \int_0^t F(\tau)\dot{x}^2_\alpha
(\tau,\alpha^*)d\tau + \int_0^T (-w''(x (\tau,\alpha^*)))\lambda_1
e^{-r\tau}x^2_\alpha (\tau,\alpha^*)d\tau \geq 0.$$ Taking into
account that $F(t)\dot{x}_\alpha (t,\alpha^*)x_\alpha
(t,\alpha^*)=\frac{d}{dt}(x^2_\alpha (t,\alpha^*))F(t)/2$, we
establish that the function $x^2_\alpha (t,\alpha^*)$ is
increasing. In view of the initial conditions, in a small interval
$(0,\varepsilon)$ we have $x_\alpha (t,\alpha^*)>0$. This
inequality holds for all $t\in (0,T)$, for otherwise there would
exist $\bar{t}\in (\varepsilon ,T)$ for which $x_\alpha
(\bar{t},\alpha^*)=0$. The latter would lead to the contradiction
$0<x^2_\alpha (\varepsilon/2,\alpha^*)\leq x^2_\alpha
(\bar{t},\alpha^*)=0$. Taking into account that $w''(x)\leq 0$, we
obtain from \eqref{eq:diffalpha1alpha2} that $\alpha_1=\alpha_2$,
i.e. the two extremals coincide.
\endprf

We augment the above results by investigating the dependence of
the final location $X(T) := x(T;T,x_0)$ of the agent on the length
of the time horizon $T$. Since $X(T)< x_1$, $\forall T >0$, we
have $l := \sup_{T>0}X(T) \leq x_1$.
\begin{proposition} Under the assumptions of Proposition \ref{thm:SingleExtremal}, we have the following classification. If $\rho \geq r$ or $\rho \in (0,r(1-\theta)]$, then $l=x_1$. If $\rho \in (r(1-\theta),r)$, it is possible to have $l < x_1$ for appropriate values of the parameters of the problem.
\label{thm:ReachingThePeak}\end{proposition} \proof Assume that
$l<x_1$. For $\rho \geq r$, we have $e^{(r-\rho)t}\leq 1$ and so
$$X(T)=x_0+ \int_0^T \frac{ \int_\tau^T \lambda_1(T) w'(x(s))e^{-rs}ds}{2 (\xi \lambda_1(T) e^{-r\tau}+
\eta e^{-\rho \tau})} d\tau \geq   x_0+ \frac{ \lambda_1(T)
w'(l)}{2 (\xi \lambda_1(T)+ \eta )} \int_0^T e^{r \tau}
\left(\int_\tau^Te^{-rs}ds \right)d\tau,$$ where the integral
evaluates to $\frac{1}{r}\left[ T-\frac{1}{r}+\frac{e^{-rT}}{r}
\right]$.

According to the results from section \ref{sec:termasset}, the
expression $\frac{ \lambda_1(T)}{2 (\xi \lambda_1(T)+ \eta )}$
does not tend to zero as $T \rightarrow \infty$, so $\lim_{T
\rightarrow \infty}X(T)=\infty$, which contradicts the fact that
$X(T)$ is bounded.

For $\rho < r$, we study three cases according to the behaviour of
$\lambda_1(T)$:
\begin{enumerate}
  \item[i)] $\rho \in (r(1-\theta),r)$. In this case $\lim_{T
\rightarrow \infty}\lambda_1(T) = const$,
  \item[ii)] $\rho \in (0,r(1-\theta))$. In this case $\lim_{T
\rightarrow \infty}\lambda_1(T) = \infty$.
  \item[iii)] $\rho = r(1-\theta)$. In this case $C_1 (1+T)^\theta \leq
  \lambda_1(T)$.
\end{enumerate}
For case i) we have $$X(T) \leq   x_0+ \frac{ \lambda_1(T)
w'(x_0)}{2 \eta} \int_0^T e^{\rho \tau} \left(\int_\tau^Te^{-rs}ds
\right)d\tau \leq x_0 + \frac{ \lambda_1(T) w'(x_0)}{2
\rho(r-\rho)\eta},$$ after taking into account that the integral
evaluates to
$\frac{1}{r}\left[\frac{1}{r-\rho}-\frac{r}{(r-\rho)\rho}e^{-(r-\rho)T}+\frac{e^{-rT}}{\rho}\right]$.
It is clear that, for instance, for large values of $\eta$ this
upper bound on $X(T)$ can be strictly smaller than $x_1$.

In case ii), defining $A := r(1-\theta)-\rho > 0$, we have from
section \ref{sec:termasset} $$C_1 e^{AT}\leq \lambda_1(T) \leq C_2
e^{AT} \textrm{ with } C_i>0,~i=1,2.$$ Consequently, assuming that
$l<x_1$, we have \beq X(T) \geq x_0+\frac{C_1 w'(l)}{2r}\int_0^T
e^{AT} \frac{e^{-r\tau}-e^{-rT}}{\zeta e^{AT}e^{-r\tau}+\eta
e^{-\rho \tau}}d\tau ,\label{eq:XTlowbound}\eeq where $\zeta :=
\xi C_2 >0 $. Denote the integral in \eqref{eq:XTlowbound} by $I$.
We have $$I = \int_0^T \frac{1-e^{r\tau}e^{-rT}}{\zeta +\eta
\frac{e^{(A+r\theta)\tau}}{e^{AT}}}d\tau \geq
\int_0^{\frac{A}{A+r\theta}T} \frac{1-e^{r\tau}e^{-rT}}{\zeta
+\eta \frac{e^{(A+r\theta)\tau}}{e^{AT}}}d\tau .$$ Since
$e^{-AT}e^{(A+r\theta)\tau}\leq 1$ for $\tau \in [0,
AT/(A+r\theta)]$, we have from the last expression $$I \geq
\int_0^{\frac{A}{A+r\theta}T} \frac{1-e^{r\tau}e^{-rT}}{\zeta
+\eta }d\tau  = \frac{A}{(A+r\theta)(\zeta+\eta)}T
-\frac{e^{-rT}}{\zeta + \eta}\cdot
\frac{e^{\frac{A}{A+r\theta}rT}-1}{r}.$$ The last expression tends
to infinity as $T \rightarrow \infty$, implying that $X(T)$ is
unbounded. This contradiction shows that $l = x_1$.

In case iii) the condition from section \ref{sec:termasset} is
equivalent to
$$\frac{1}{\lambda_1(T)}\leq \frac{1}{C_1 (1+T)^\theta}.$$ Assume
that $l < x_1$. Then
\begin{equation*}
\begin{split}
X(T) \geq & x_0+\frac{w'(l)\lambda_1(T)}{2r}\int_0^T
\frac{1-e^{-rT}e^{r\tau}}{\xi \lambda_1(T)+\eta
e^{r\theta\tau}}d\tau \geq  x_0 + \frac{w'(l)}{2r}\int_0^T
\frac{1-e^{-rT}e^{r\tau}}{\xi + \frac{\eta e^{r\theta\tau}}{C_1
(1+T)^\theta}}d\tau = \\ & x_0 + \frac{w'(l)}{2r\eta}C_1
(1+T)^\theta \int_0^T \frac{1-e^{-rT}e^{r\tau}}{\xi \frac{C_1
(1+T)^\theta}{\eta} +  e^{r\theta\tau}}d\tau .
\end{split}
\end{equation*}
Set $B=B(T):= \xi \frac{C_1 (1+T)^\theta}{\eta}$ and introduce the
change of variables $\mu=e^{r\theta\tau}$ in the last expression
to obtain $$X(T) \geq x_0+const (1+T)^\theta \int_1^{e^{r\theta
T}} \frac{1-e^{-rT}\mu^{\frac{1}{\theta}}}{r\theta \mu
(\mu+B)}d\tau.$$ This requires us to study the behaviour of two
expressions.

First, we have $$ (1+T)^\theta \int_1^{e^{r\theta T}} \frac{d
\mu}{\mu(\mu + B)} = (1+T)^\theta \frac{1}{B}\left[ \ln\left(
\frac{1}{1+\frac{const (1+T)^\theta}{e^{r\theta T}}} \right) + \ln
\left( 1+const (1+T)^\theta \right) \right].$$ When $T \rightarrow
\infty$, the first logarithm tends to zero and the second one
tends to infinity, i.e. the whole expression tends to infinity.

Second, note that
$$0 \leq (1+T)^\theta e^{-rT} \int_1^{e^{r\theta T}} \frac{\mu^{\frac{1}{\theta}}}{\mu(\mu + B)}d
\mu \leq \frac{(1+T)^\theta}{ e^{rT}} \int_1^{e^{r\theta T}}
\frac{\mu^{\frac{1}{\theta}}}{\mu^2}d \mu =
\frac{\theta}{1-\theta}\left[ \frac{(1+T)^\theta}{e^{r\theta
T}}-\frac{(1+T)^\theta}{e^{r T}} \right].$$ The last expression
tends to zero as $T \rightarrow \infty$.

Combining the above results, we obtain $X(T) \rightarrow \infty$,
which contradicts the fact that $X(T)$ is bounded. Thus, in this
case $l = x_1$. \endprf

\end{document}